\documentclass[a4paper,12pt]{article}

\usepackage[plainpages=false]{hyperref}
\usepackage[english,ngerman]{babel}
\usepackage{abstract} 
\usepackage[latin1]{inputenc}
\usepackage{verbatim}
\usepackage[arrow, matrix, curve]{xy}
\usepackage{amsmath,amssymb,amsfonts,amsthm,amsopn}
\usepackage{nicefrac}
\usepackage{fancyhdr}

\title{Heights and totally real numbers}
\author{Lukas Pottmeyer \vspace{0.4cm}}
\date{}

\DeclareMathOperator{\supp}{supp}
\DeclareMathOperator{\Gal}{Gal}
\DeclareMathOperator{\PrePer}{PrePer}

\newtheorem{Theorem}[]{Theorem}[section]
\newtheorem{maintheorem}{Theorem}

\newtheorem{corollary}[Theorem]{Corollary}
\newtheorem{Proposition}[Theorem]{Proposition}
\newtheorem{Lemma}[Theorem]{Lemma}
\newtheorem*{Definition}{Definition}
\newtheorem{Ohne}[Theorem]{}
\newtheorem{Example}[]{Example}
\newtheorem{Question}[]{Question}
\newtheorem{Remark}[Theorem]{Remark}
\newtheorem{Facts}[Theorem]{Facts}

\makeindex

\begin{document}
 
 \maketitle
 
 \selectlanguage{english}
 \begin{abstract} In 1973 Schinzel proved in \cite{Sch73} that the standard logarithmic height $h$ on the maximal totally real field extension of the rationals is either zero or bounded from below by a positive constant. In this paper we study this property for canonical heights associated to rational functions and the corresponding dynamical system on the affine line.
 \end{abstract}
 
 
 \begin{section}{Introduction}\label{Introduction}
 
 We fix an algebraic closure $\overline{\mathbb{Q}}$ and we denote the maximal totally real algebraic subfield by $\mathbb{Q}^{tr}$.
Let $h$ be the standard logarithmic height on the algebraic numbers. We say that a field $F\subset \overline{\mathbb{Q}}$ has the Bogomolov property relative to $h$ if and only if $h(\alpha)$ is either zero or bounded from below by a positive constant for all $\alpha \in F$. This notation was introduced 2001 by Bombieri and Zannier in \cite{BZ01}. The name is given in analogy to the famous Bogomolov conjecture, yielding a lower bound of the N\'eron-Tate height on a certain set of algebraic points on an abelian variety (see \cite{BG}, Theorem 11.10.17). By Northcott's theorem, every number field has the Bogomolov property relative to $h$. In Table  \ref{Bogomolov h} we summerize some examples of fields $F$ of infinite degree over $\mathbb{Q}$ with the Bogomolov property relative to $h$. Let $K$ be a number field, then we denote by $K^{ab}$ the maximal abelian field extension of $K$. Furthermore, a field is called totally $p$-adic if and only if it may be embedded in a finite extension of $\mathbb{Q}_p$. This is a $p$-adic analogue of the field $\mathbb{Q}^{tr}$.

\begin{table}[ht]
	\centering
		\begin{tabular}{| c | c |}
		 \hline
			Field & Reference \\ \hline\hline
			$\mathbb{Q}^{tr}$ & Schinzel \cite{Sch73} \\ \hline
			finite extensions of $K^{ab}$ & Amoroso, Zannier \cite{AZ00} \\ \hline
			totally $p$-adic fields & Bombieri, Zannier \cite{BZ01} \\ \hline
			$\mathbb{Q}(E_{\rm tor})$, $E/\mathbb{Q}$ elliptic curve & Habegger \cite{Ha11} \\ \hline
		\end{tabular}
	\caption{Fields with the Bogomolov property relative to $h$}
	\label{Bogomolov h}
\end{table}

We want to study the behavior of canonical heights associated to rational functions on $\mathbb{Q}^{tr}$, and hence variations of Schinzel's result. For convenience we will state it again as a theorem.
 
 \begin{Theorem}[Schinzel]\label{Schinzel}
The field $\mathbb{Q}^{tr}$ has the Bogomolov property relative to $h$.
\end{Theorem}

Schinzel did not use the notation ``Bogomolov property'' in his paper. As we have mentioned above, this notation was introduced nearly 30 years after Schinzel's result. The proof of Schinzel gives the sharp lower bound $\frac{1}{2}\log(\frac{1+\sqrt{5}}{2})$. If one is interested only in a quantitative result, one can use Bilu's equidistribution theorem (see \cite{BG}, Theorem 4.3.1) as follows:

Assume there is a sequence $\{\alpha_n \}_{n \in \mathbb{N}}$ in $\mathbb{Q}^{tr}\setminus \{-1,0,1\}$, such that the height of these points tends to zero. Then, by Bilu's equidistribution theorem, the equidistributed probability measures on the set of conjugates of the $\alpha_i$ converge weakly to the probability measure on the unit circle. But the support of each such Galois measure lies in the real line. Hence they cannot cover the unit circle, which leads to a contradiction.

See also \cite{HS93} for a very short proof of Schinzel's original result. Roughly speaking, the field $\mathbb{Q}^{tr}$ seems to be arithmetically ``easy". In this paper we will give a complete classification of rational functions defined over the algebraic numbers such that $\mathbb{Q}^{tr}$ has the Bogomolov property relative to the canonical height coming from this rational map. For a rational function $f$ we denote by $\PrePer(f)$ the set of preperiodic points of $f$; i.e. points with finite forward orbit. Our result reads as follows.

\begin{maintheorem}\label{rational maps}
Let $f\in \overline{\mathbb{Q}}(x)$ be a rational function of degree at least two. Then the following statements are equivalent:
\begin{itemize}
	\item[i)] $\mathbb{Q}^{tr}$ has the Bogomolov property relative to $\widehat{h}_f$.
	\item[ii)] There is a $\sigma \in \Gal(\overline{\mathbb{Q}}/\mathbb{Q})$ such that the Julia set of $\sigma(f)$ is not contained in $\mathbb{R}$.
	\item[iii)] The set $\PrePer(f)\cap \mathbb{Q}^{tr}$ is finite.
\end{itemize}
\end{maintheorem}

\vspace{0.3cm}

Statements $i)$ and $iii)$ will follow directly from $ii)$ by the equidistribution theorem of Yuan which is a far-reaching generalization of Bilu's theorem.

In Section \ref{Heights and dynamics} we give a short introduction to canonical dynamical heights and state a few facts on Julia sets that will be needed in the proof of Theorem \ref{rational maps}. Section \ref{Example and Lemma} contains a proof of a partial result of our main theorem. This result is completely covered by Theorem \ref{rational maps}, but the proof is very simple and shows the strategy for proving Theorem \ref{rational maps} very clearly. In Section \ref{Proof} we prove the main theorem and give an additional equivalence in the case of a polynomial. One class of polynomials with real Julia set are Chebyshev polynomials. We will briefly study these polynomials in Section \ref{Chebyshev polynomials}.

\vspace{0.3cm}

\textit{Acknowledgment:} Most of this work was done during a research stay at the Institute for computational and experimental research in mathematics (ICERM) in Providence. I am thankful for their hospitality during February and March 2012, and for the support of the DFG-Graduiertenkolleg GRK 1692. Furthermore I would like to thank Fabrizio Barroero, Paul Fili, Walter Gubler, Khoa Nguyen, Adam Towsley, Tom Tucker and Umberto Zannier for lively discussions and very helpful remarks and suggestions. 
Moreover, I thank the referee for the very careful reading of a first version of this paper.
\end{section}

\begin{section}{Heights and dynamical Systems}\label{Heights and dynamics}

Canonical heights associated to rational functions defined over the algebraic numbers can be defined using the next theorem due to Call and Silverman.
  
 \begin{Theorem}\label{canonical height}
 Let $f\in\overline{\mathbb{Q}}(x)$ be a rational function of degree greater one. There is a unique height function $\widehat{h}_f$, called the canonical height related to $f$, such that for all $\alpha\in\overline{\mathbb{Q}}$ we have
 $$\begin{array}{lcr}
 i) \hspace{0.2cm} \widehat{h}_f (f(\alpha))=\deg(f)\widehat{h}_f (\alpha) & \text{ and } & ii) \hspace{0.2cm} \widehat{h}_f = h + O(1) .
 \end{array}$$
The canonical height $\widehat{h}_f$ vanishes precisely on the set $\PrePer(f)$.
 \end{Theorem}
 
See \cite{Si07}, Chapter 3.4, for a proof and additional information including the following two properties which we will use frequently. 

\begin{Proposition}\label{Northcott} With the notation from Theorem \ref{canonical height} we have
\begin{itemize}
\item[a)] $\widehat{h}_f (\alpha)=0$ $\Leftrightarrow$ $\alpha\in\PrePer(f)$,
\item[b)] $\vert \{\alpha \vert \widehat{h}_f (\alpha) \leq A, \deg(\alpha)\leq B \}\vert<\infty$ for all $A,B \in \mathbb{R}$.
 \end{itemize}
 \end{Proposition}

Proposition \ref{Northcott} $b)$ is commonly known as Northcott's theorem.

\begin{Definition}\rm
Let $f\in\overline{\mathbb{Q}}(x)$ be a rational function of degree at least $2$. We say that a field $F\subset\overline{\mathbb{Q}}$ has the Bogomolov property relative to $\widehat{h}_f$, if and only if there exists a positive constant $c$ such that $\widehat{h}_f (\alpha)\geq c$ for all $\alpha\in F\setminus\PrePer(f)$.
\end{Definition}
 
 Notice that the standard height $h$ fulfills $h(\alpha^d)=dh(\alpha)$, for all $\alpha\in\overline{\mathbb{Q}}$ and all $d\in\mathbb{N}$. Hence we find $$\widehat{h}_{x^d}=h \text{ for all } d\geq 2 \quad,$$  and for $f=x^d$ the above definition coincides with the definition given by Bombieri and Zannier in \cite{BZ01}.
 
 In fact, we work with rational functions on the Riemann sphere which we identify with $\mathbb{C}\cup\{\infty\}$.  
 On the Riemann sphere, we will always use the complex topology which is induced by the chordal metric $\rho$. Let $z\in\mathbb{C}$, recall that the chordal metric is given by $$\rho(z,z'):= \begin{cases}\frac{\vert z-z'\vert}{\sqrt{1+\vert z\vert^2}\sqrt{1+\vert z' \vert^2}}& \text{ if } z'\neq \infty \\ \frac{1}{\sqrt{1+\vert z \vert^2}}& \text{ if } z' = \infty \end{cases}\quad.$$ The Julia set of such a map $f$ is the set of points where $f$ acts "chaotically".

\begin{Definition}\rm
Let $f$ be a self map of the Riemann sphere. The Fatou set $F(f)$ of $f$ is the maximal open subset of the Riemann sphere, satisfying the condition: For all $\alpha \in F(f)$ and all $\varepsilon > 0$ there exists a $\delta>0$ such that $$\rho(\alpha , \beta) < \delta \Rightarrow \rho(f^n (\alpha), f^n (\beta)) < \varepsilon$$ for all $n\in\mathbb{N}$. The Julia set $J(f)$ of $f$ is the complement of $F(f)$.
\end{Definition}
 
 In addition to the canonical height associated with a rational function $f$ of degree $\geq2$ there exists a $f$-invariant canonical probability measure $\mu_f$ which is supported on the Julia set of $f$ (see \cite{FLM83}).
 
\begin{Ohne}\rm Although we do not need the theory of polarized algebraic dynamical systems for our main theorem, we will briefly recall the basic definitions. For detailed information we refer to \cite{Yu12} and the references therein. 

Let $K$ be a number field, and let $X$ be a smooth projective variety of dimension $n$ with a morphism $f: X \rightarrow X$, both defined over $K$. Moreover, let $L$ be an ample line bundle on $X$. The triple $(X,L,f)$ is called (polarized) algebraic dynamical system if we have $f^{*}L\cong L^{\otimes q}$, for $q \geq 2$. We need to fix a line bundle $L$ to associate a canonical height and a canonical measure to the algebraic dynamical system. The canonical height $\widehat{h}_{X,L,f}$ for $(X,L,f)$ is uniquely determined by the properties given in Theorem \ref{canonical height}. Namely, 
$$\begin{array}{lcr}
 \widehat{h}_{X,L,f} (f(P))=q\widehat{h}_{X,L,f} (P) \hspace{0.2cm} \forall P\in X(\overline{K}) & \text{ and } & \hspace{0.2cm} \widehat{h}_{X,L,f} = h_L + O(1) ,
 \end{array}$$
where $h_L$ is any Weil height on $X$ (see \cite{CS93}). 

For a fixed non-archimedean $v\in M_K$ we write $\mathbb{C}_v$ to denote the completion of an algebraic closure of $K_v$. This is a complete and algebraically closed field (see \cite{BGR}, Proposition 3.4.3). We consider $(X,L,f)$ as an algebraic dynamical system defined over $\mathbb{C}_v$ and denote by $X_{v}^{an}$ the associated Berkovich space to $X/\mathbb{C}_v$. For the theory of Berkovich spaces we refer to \cite{Ber}. As in the special case above, there exists a $v$-adic canonical $f$-invariant measure $\mu_{f,v}$ on $X_{v}^{an}$ associated to $(X,L,f)$. This is also true for archimedean $v\in M_K$, where we set $X_{v}^{an}:=X(\mathbb{C})$ regarded as a complex manifold. 

The canonical height associated to the algebraic dynamical system $(\mathbb{P}_{K}^{1},\mathcal{O}(1),f)$, $f\in K(x)$ of degree $\geq 2$, is the function $\widehat{h}_f$ from Theorem \ref{canonical height}. The map $f$ extends uniquely to a continuous function on $\mathbb{P}^{1}(\mathbb{C}_v)$, $v \in M_K$. We define the Berkovich Julia set $J^{\mathcal{B}}_v (f)$ of $f$, to be the support of the canonical measure $\mu_{f,v}$. 
\end{Ohne}

Let $P\in X(\overline{\mathbb{Q}})$ be arbitrary and let $\delta_P$ be the Dirac measure at $P$. We denote the set $\{ \sigma(P) \vert \sigma \in \Gal(\overline{\mathbb{Q}}/K) \}$ of $K$-Galois conjugates of $P$ by $G_K (P)$, and define the probability measure $$\overline{\delta_P}:=\vert G_K (P) \vert^{-1} \sum_{P' \in G_K (P)} \delta_{P'}\quad .$$ Now we can formulate Yuan's equidistribution theorem (see \cite{Yu08}, Theorem 3.7).

\begin{Theorem}[Yuan]\label{Yuan}
Let $(X,L,f)$ be a polarized algebraic dynamical system defined over the number field $K$, and let $\{ P_i \}_{i \in\mathbb{N}}$ be a sequence of pairwise distinct points in $X(\overline{\mathbb{Q}})$ such that
\begin{itemize}
\item[i)] $\widehat{h}_f (P_i ) \rightarrow 0$, as $i \rightarrow \infty$,
\item[ii)] every infinite subsequence of $\{ P_i \}_{i \in\mathbb{N}}$ is Zariski dense in $X$.
\end{itemize}
For any $v\in M_K$ the measures $\overline{\delta_i}:=\overline{\delta_{P_i}}$ converge weakly to $\mu_{f,v}$. This means that for every continuous function $\varphi: X_{v}^{an} \rightarrow \mathbb{C}$ we have
$$\int_{X_{v}^{an}} \varphi(x) \overline{\delta_i}=\vert G_K (P_i )\vert^{-1} \sum_{P'_i \in G_K (P_i )} \varphi(P_i ) \rightarrow \int_{X_{v}^{an}}\varphi(x)\mu_{f,v} \quad ,$$
as $i \rightarrow \infty$.
\end{Theorem}

Of course, the second requirement on the sequence $\{ P_i \}_{i \in\mathbb{N}}$ in the above theorem is always true if $X=\mathbb{P}^{1}$.

We want to study canonical heights $\widehat{h}_f$ on the field $\mathbb{Q}^{tr}$. This was our main motivation for the next theorem and a first version only covers Corollary \ref{totreal Bogomolov}. Paul Fili pointed out that the same proof applies in a more general setting (see \cite{FM12}).
 
\begin{Theorem}\label{v-adic Bogomolov}
Let $f \in \overline{\mathbb{Q}}(x)$ be a rational function of degree $\geq2$, and let $K$ be a number field with valuation $v\in M_K$ such that the Berkovich Julia set $J^{\mathcal{B}}_{v} (f)$ is not contained in the closure of $(\mathbb{P}^{1})^{an}_{v} (K)$. If $L/K$ is a Galois extension lying in $K_v$, then $L$ has the Bogomolov property relative to $\widehat{h}_f$. Furthermore, there are only finitely many preperiodic points of $f$ in $L$.
\end{Theorem}
 
\textbf{Proof:} Let $F$ be a number field such that $f\in F(x)$ and $K\subseteq F$. Assume there is a sequence $\{\alpha_i \}_{i\in{\mathbb{N}}}$ in $L$ of pairwise distinct elements satisfying $\widehat{h}_f (\alpha_i ) \rightarrow 0$ for $i\rightarrow \infty$. Denote by $\overline{\delta_i}$ the equidistributed probability measures on the set $G_F (\alpha_i )$. The support of $\overline{\delta_i}$ lies in $K_v$ for all $i \in \mathbb{N}$, as $L/K$ was assumed to be Galois. Notice that the choice of $F$ implies $G_F (\alpha_i ) \subseteq G_K (\alpha_i )$, for all $i\in\mathbb{N}$.

By assumption, there exists an $\alpha \in J^{\mathcal{B}}_{v}(f) = \supp(\mu_{f,v} )$ which is not contained in the closure of $K_v$ in $(\mathbb{P}^{1})^{an}_{v}$. As $(\mathbb{P}^{1})_{v}^{an}$ is a Hausdorff space, there is an open neighborhood $U$ of $\alpha$ such that $U\cap K_v = \emptyset$. By Theorem \ref{Yuan}, the measures $\overline{\delta_i}$ converge weakly to $\mu_{f,v}$ and hence $$0=\lim_{i\rightarrow \infty} \overline{\delta_i}(U) = \mu_{f,v}(U)\neq 0 \quad .$$ This is a contradiction, and hence there cannot exist such a sequence $\{\alpha_i\}_{i\in\mathbb{N}}$. \hfill $\square$ 

\vspace{0.3cm}

The case $K=\mathbb{Q}$, $L=\mathbb{Q}^{tr}$ and $v=\infty$ yields the following result.

\begin{corollary}\label{totreal Bogomolov}
Let $f \in \overline{\mathbb{Q}}(x)$ be a rational function of degree $\geq 2$ such that the Julia
set $J(f)$ of $f$ is not contained in the real line. Then $\mathbb{Q}^{tr}$ has the Bogomolov property
relative to $\widehat{h}_f$. Furthermore, there are only finitely many preperiodic points of $f$ in $\mathbb{Q}^{tr}$.
\end{corollary}
 
 Notice again that $h=\widehat{h}_{x^2}$ and the Julia set of the map $x^2$ is indeed the unit circle. Hence Schinzel's result (in a non-effective version) is a special case of Corollary \ref{totreal Bogomolov}.

\begin{Remark}\rm
Corollary \ref{totreal Bogomolov} also includes a special case of a theorem of Zhang (\cite{Zh98}, Corollary 2). Let $E$ be an elliptic curve defined over a number field $K$ with N\'eron-Tate height $\widehat{h}_E$. Then there exists a rational map $f\in K(x)$, called Latt\`es map, such that $\widehat{h}_E (P) = \frac{1}{2}\widehat{h}_f (x(P))$ for all $P \in E(\overline{\mathbb{Q}})$. Further, the Julia set of $f$ is the complete Riemann sphere. Hence Theorem \ref{totreal Bogomolov} tells us that there is a positive constant $c$ such that $\widehat{h}_E (P) \geq c$ for all non-torsion points $P \in E (\mathbb{Q}^{tr})$ and there are only finitely many torsion points in $E (\mathbb{Q}^{tr})$. Notice that there is an effective constant $c$ in the case where $K$ is totally real (see \cite{BP05}, Theorem 17).
\end{Remark}

In the non-archimedean case, Theorem \ref{v-adic Bogomolov} gives a dynamical version of the result of Bombieri and Zannier stated in the introduction (\cite{BZ01}, Theorem 2). For details on this case we refer to \cite{FM12}. 
 
We will collect some important facts on Julia sets of rational maps.
 
 \begin{Facts}\label{Julia properties}
 Let $f \in\mathbb{C}(x)$ be a rational function of degree at least two. Then we have
\begin{itemize}
	\item[a)] $J(f)$ is not empty,
	\item[b)] $J(f)$ is completely invariant, i.e. $f(J(f))=f^{-1}(J(f))=J(f)$,
	\item[c)] there are no isolated points in $J(f)$,
	\item[d)] $J(f)$ is the closure of the repelling periodic points of $f$,
	\item[e)] The backward orbit of any $\alpha \in J(f)$ (i.e. $\cup_{n\in\mathbb{N}}\{\beta\in\overline{\mathbb{Q}} \vert f^{n}(\beta)=\alpha\}$) is dense in $J(f)$.
\end{itemize}
 \end{Facts}
 
 \textbf{Proof:} For proofs of these statements we refer to \cite{Be}, Theorem 4.2.1, Theorem 3.2.4, Theorem 5.7.1, Theorem 6.9.2 and Theorem 4.2.7 $ii)$. \hfill $\square$

\end{section}

\begin{section}{A first example}\label{Example and Lemma}

 A natural question arising from Corollary \ref{totreal Bogomolov} is the following. Does $\mathbb{Q}^{tr}$ have the Bogomolov property relative to $\widehat{h}_f$ even if $J(f)$ is contained in the real line? Before proving our main theorem which gives a complete answer to this question, we will give a counterexample. Therefore we need a map whose Julia set is contained in the real line. If we restrict to the polynomial case, we have the following lemma.

\begin{Lemma}\label{real polynomial}
Let $f\in\mathbb{C}[x]$ be a polynomial of degree at least $2$. Then we have $J(f)\subset\mathbb{R}$ if and only if all preperiodic points of $f$ are real.
\end{Lemma}

\textbf{Proof:} If all preperiodic points of $f$ are real, then in particular the closure of the repelling periodic points of $f$ lies in the real line. This set is just $J(f)$ (see Fact \ref{Julia properties} $d)$). Next we assume $J(f)\subset\mathbb{R}$. As we have a polynomial, the Julia set of $f$ is the boundary of the compact set $$\{ y \in \mathbb{C} \vert ~ \vert f^n (y) \vert \nrightarrow \infty \text{, as } n \rightarrow \infty \} \quad .$$
See \cite{Mi}, Lemma 9.4. This set is called the filled Julia set of $f$. For a polynomial it follows from the definitions of the Julia set and the filled Julia set that $\infty$ is in neither of both sets. Hence, both sets are bounded. By assumption, $J(f)$ is a closed subset of a closed interval $I$. The only bounded subset of the Riemann sphere with such a boundary is the set itself, proving that $J(f)$ coincides with the filled Julia set. Of course all preperiodic points of $f$ are contained in the filled Julia set, proving the lemma. \hfill $\square$

\begin{Example}\label{classification quadratic maps} \rm
The Julia set of $f_c (x)=x^2 -c \in \mathbb{C}[x]$ is contained in the real line if and only if $c\in \mathbb{R}$ and $c \geq 2$. To prove this claim we assume first that $c$ is either non-real or $<2$. Then $\sqrt{-\nicefrac{1}{2}\sqrt{4c+1}-\nicefrac{1}{2}+c}$ is a non-real preperiodic point of $f_c$, and from Lemma \ref{real polynomial} we know $J(f_c )\not\subset\mathbb{R}$. If $c\geq 2$, then the pre-image of the interval $[-c,c]$ is contained in this interval. Induction yields $f_{c}^{-n}([-c,c])\subseteq[-c,c]$ for all $n\in\mathbb{N}$. Moreover, the fixed point $\nicefrac{1}{2}\sqrt{4c+1}+\nicefrac{1}{2}$ of $f_c$ is repelling and lies in the interval $[-c,c]$. By Facts \ref{Julia properties} $d)$ and $e)$, we conclude $J(f_c )\subseteq[-c,c]$.
\end{Example}
 
 \begin{Proposition}\label{noBogomolov}
Let $c$ be a rational with $c\geq2$. Then $\mathbb{Q}^{tr}$ does not have the Bogomolov property relative to $\widehat{h}_{f_c}$.
 \end{Proposition}
 
\textbf{Proof:} Take an $\epsilon \in (-c,c)\cap \mathbb{Q}$ such that $\epsilon$ is not a preperiodic point of $f_c$. This is possible by Northcott's theorem (see Proposition \ref{Northcott}). As seen in the example, the set $f_{c}^{-n}(\epsilon)$ is contained in $\mathbb{R}$ for all $n\in\mathbb{N}$. Moreover, the sets $f_{c}^{-n}(\epsilon)$ are Galois invariant, since $c$ and $\epsilon$ were chosen to be rational numbers. Hence, $f_{c}^{-n}(\epsilon)$ is actually contained in $\mathbb{Q}^{tr}$.
For all $n$ we take an arbitrary $\gamma_n$ in $f_{c}^{-n}(\epsilon)$ and get a sequence $\{\gamma_n\}_{n\in \mathbb{N}}$ in $\mathbb{Q}^{tr}$, with $\widehat{h}_{f_c}(\gamma_n)=\frac{1}{2^n}\widehat{h}_{f_c}(\epsilon)$. This tends to zero, proving the claim. Notice that the canonical height of these $\gamma_n$ is positive, as $\epsilon$ (and thus $\gamma_n$) are not preperiodic. \hfill $\square$

\end{section}

\begin{section}{Proof of the main result}\label{Proof}

We already stated in the Introduction that Yuan's equidistribution theorem will play an essential role in the proof of our main theorem. Another important ingredient is the following theorem of Eremenko and van Strien (see \cite{EvS11}, Theorem 2 and the discussion afterwards).

\begin{Theorem}\label{Eremenko vanStrien}
Let $f\in\mathbb{C}(x)$ be an rational function of degree at least $2$. If $J(f)\subset\mathbb{R}$, then there exist finitely many intervals $I_1 , \dots ,I_r$ such that $f^{-1}(I_1 \cup \dots \cup I_r) \subseteq I_1 \cup \dots \cup I_r$.
\end{Theorem}

\begin{Remark}\label{real coefficients}\rm
Let $f$ be a rational function defined over $\mathbb{C}$ such that $\vert f(\mathbb{R})\cap\mathbb{R})\vert =\infty $. Write $f=\frac{p}{q}$, with polynomials $p, q\in\mathbb{C}[x]$ and $p$ monic. Then $f$ is actually defined over $\mathbb{R}$. To prove this claim we denote by $\overline{q}$ the coefficient-wise complex conjugation of $q$. Then $f= \frac{p\overline{q}}{q\overline{q}}=\frac{p_1 }{r} + i\frac{p_2 }{r}$, with real polynomials $p_1,p_2,r$. By assumption, $p_2$ must be zero and hence $f$ is defined over $\mathbb{R}$. In particular, $f$ is defined over $\mathbb{R}$ whenever $J(f)\subseteq\mathbb{R}$ (see Fact \ref{Julia properties} $b)$).
\end{Remark}

Every $\sigma$ in $\Gal(\overline{\mathbb{Q}}/\mathbb{Q})$ extends to a unique endomorphism of $\overline{\mathbb{Q}}(x)$ with $\sigma(x)=x$. So we can define the rational map $\sigma(f)$ for all $f\in\overline{\mathbb{Q}}(x)$ and all $\sigma \in \Gal(\overline{\mathbb{Q}}/\mathbb{Q})$ . 

\begin{Lemma}\label{galoisheights}
Let $f\in \overline{\mathbb{Q}}(x)$ be a rational map of degree $>1$. Then we have $\widehat{h}_f = \widehat{h}_{\sigma(f)} \circ \sigma$.
\end{Lemma}

\textbf{Proof:} This follows directly from the definition and the trivial facts $\deg(f) = \deg(\sigma(f))$ and $\sigma(f(\alpha))=\sigma(f)(\sigma(\alpha))$ for all $\alpha \in \overline{\mathbb{Q}}$. \hfill $\square$

\vspace{0.3cm}

Now we are prepared to prove our main theorem which we will state again for the reader's convenience.

\newtheorem*{rationalmaps}{Theorem \ref{rational maps}}
\begin{rationalmaps}
As usual let $f\in \overline{\mathbb{Q}}(x)$ be a rational map of degree at least two. Then the following statements are equivalent:
\begin{itemize}
	\item[i)] $\mathbb{Q}^{tr}$ has the Bogomolov property relative to $\widehat{h}_f$.
	\item[ii)] There is a $\sigma \in \Gal(\overline{\mathbb{Q}}/\mathbb{Q})$, such that the Julia set $J(\sigma(f))$ is not contained in $\mathbb{R}$.
	\item[iii)] The set $\PrePer(f)\cap \mathbb{Q}^{tr}$ is finite.
\end{itemize}
\end{rationalmaps}

\textbf{Proof:} Notice again that $J(f)$ cannot be empty, see Fact \ref{Julia properties} $a)$. By Corollary \ref{totreal Bogomolov}, we will conclude easily that $ii)$ yields $i)$ and $iii)$. Assume there is a $\sigma \in \Gal(\overline{\mathbb{Q}}/\mathbb{Q})$ such that $J(\sigma(f))$ is not contained in the real line. 
Then Corollary \ref{totreal Bogomolov} implies that $\mathbb{Q}^{tr}$ has the Bogomolov property relative to $\widehat{h}_{\sigma(f)}$. By Lemma \ref{galoisheights} we know that $\sigma^{-1}(\mathbb{Q}^{tr})=\mathbb{Q}^{tr}$ has the Bogomolov property relative to $\widehat{h}_f$ as well, which yields $i)$. Notice that $\mathbb{Q}^{tr}$ is a Galois extension of $\mathbb{Q}$.
Moreover, we have $\vert \PrePer(f)\cap \mathbb{Q}^{tr} \vert=\vert \sigma^{-1}(\PrePer(\sigma(f))\cap \mathbb{Q}^{tr}) \vert$. This is a finite set by Corollary \ref{totreal Bogomolov}, proving $iii)$.

\vspace{0.2cm}

The implication $iii) \Rightarrow ii)$ is not hard either. Assume $J(\sigma(f))$ is contained in the real line for all $\sigma \in \Gal(\overline{\mathbb{Q}}/\mathbb{Q})$. Using Facts \ref{Julia properties}, we see that $J(f)$ contains the infinite set of repelling periodic points of $f$. For all maps $\sigma(f)$, $\sigma \in \Gal(\overline{\mathbb{Q}}/\mathbb{Q})$, there are only finitely many non-repelling periodic points (see \cite{Be}, 9.6). Hence there are infinitely many points $\alpha\in \overline{\mathbb{Q}}$ such that $\sigma(\alpha)$ is a repelling periodic point of $\sigma(f)$, for all $\sigma \in \Gal(\overline{\mathbb{Q}}/\mathbb{Q})$. It follows from our assumption that all these $\alpha$ are totally real. In particular we get $\vert \PrePer(f)\cap \mathbb{Q}^{tr} \vert = \infty$.

\vspace{0.2cm}

Finally, we prove that $i)$ implies $ii)$. Assume again that $J(\sigma(f))$ is contained in the real line for all $\sigma\in\Gal(\overline{\mathbb{Q}}/\mathbb{Q})$. As in Remark \ref{real coefficients} we find that $f\in K(x)$ for a totally real number field $K$. Let $\sigma_1 , \dots , \sigma_d$ be a complete set of embeddings of $K$ into $\overline{\mathbb{Q}}$. By Theorem \ref{Eremenko vanStrien} we have for each $\sigma_i (f)$ a finite set of intervals such that all backward orbits of these intervals again lie in this finite set of intervals. 
Thus, for all $\sigma_i$ we can choose a real interval $(a_i, b_i)$ such that for all $c\in(a_i, b_i)$ every backward orbit is contained in the real line. For all $\sigma_i$ take a $c_{i}\in(a_i, b_i)\cap\mathbb{Q}$ and choose a global $\varepsilon >0$ such that $( c_i -\varepsilon , c_i + \varepsilon )\subset (a_i, b_i)$ for all $1\leq i \leq d$. All the $\sigma_i$ give rise to non-equivalent absolute values on $K$. By the approximation theorem of Artin and Whaples (see \cite{La}, Chapter II, 1), there exists a $c \in K$, such that $\vert \sigma_i (c-c_i )\vert=\vert \sigma_i (c) - c_i \vert < \varepsilon$. This implies that $\sigma_i (c)$ lies in the interval $(a_i, b_i)$ for all $\sigma_i$. For this conclusion we used the fact that $K$, and hence $c$, is totally real. There are infinitely many points $c$ with this property in $K$, but as a number field $K$ contains only finitely many preperiodic points of $f$ (by Northcott's theorem, Proposition \ref{Northcott} $b)$). Thus we can assume that $c$ is no preperiodic point of $f$.

For every $\gamma$ with $f^n (\gamma)=c$ we have $\sigma(f)^n (\sigma(\gamma)) = \sigma(c)$, $n \in \mathbb{N}$. From the choice of our intervals it follows that all conjugates of $\gamma$ are in the real line, and hence we can conclude $f^{-n}(c) \subset \mathbb{Q}^{tr}$. Now choose for all $n\in\mathbb{N}$ a $\gamma_n$ in $f^{-n}(c)$. This gives a sequence $\{\gamma_n\}$ in $\mathbb{Q}^{tr}$, such that
$$0 \neq \widehat{h}_f (\gamma_n)=\frac{1}{\deg(f)^n} \widehat{h}_f (c) \rightarrow 0 \quad.$$
Notice that we have chosen a non-preperiodic $c$. This shows that $\mathbb{Q}^{tr}$ cannot have the Bogomolov property relative to $\widehat{h}_f$. \hfill $\square$

\vspace{0.3cm}  

In the case where $f\in \overline{\mathbb{Q}}[x]$ is a polynomial we can give a further nice equivalence.

\begin{corollary}\label{equivalence polynomials}
Let $f\in \overline{\mathbb{Q}}[x]$ be a polynomial. Then the following statements are equivalent:
\begin{itemize}
	\item[i)] $\mathbb{Q}^{tr}$ does not have the Bogomolov property relative to $\widehat{h}_f$
	\item[ii)] $J(\sigma(f))\subset\mathbb{R}$ for all $\sigma\in\Gal(\overline{\mathbb{Q}}/\mathbb{Q})$
	\item[iii)] $\PrePer(f)\subset\mathbb{Q}^{tr}$
	\item[iv)] $\widehat{h}_f (\alpha) > 0$ for all $\alpha \in \overline{\mathbb{Q}}\setminus\mathbb{Q}^{tr}$ 
\end{itemize}
\end{corollary}

\textbf{Proof:} $i)$ and $ii)$ are equivalent by Theorem \ref{rational maps} and the equivalence of $iii)$ and $iv)$ is trivial by Proposition \ref{Northcott} $a)$. 
The equivalence of $ii)$ and $iii)$ follows from Lemma \ref{real polynomial} and the fact $\sigma(\PrePer(f))=\PrePer(\sigma(f))$ for all $\sigma$ in the absolute Galois group of $\mathbb{Q}$. \hfill $\square$

\begin{Remark}\rm
The Bogomolov property is in general not preserved under finite field extension. The known counterexample (see \cite{AN07} and \cite{ADZ11}) is the extension $\mathbb{Q}^{tr}(i)$ which does not have the Bogomolov property relative to the standard height $h$. We can prove this fact using dynamical methods and Theorem \ref{rational maps}. The M\"obius transformation $g(x)=\frac{x+i}{x-i}$ maps the real line onto the unit circle. Take the map $g^{-1}\circ x^2 \circ g$. By \cite{Be}, Theorem 3.1.4, we have $J(g^{-1}\circ x^2 \circ g)=g^{-1}(J(x^2))=\mathbb{R}$. The same is true for the only Galois conjugate $\frac{x-i}{x+i}$ of $g$. Furthermore, it is easy to check that we have $\widehat{h}_{g^{-1}\circ x^2 \circ g}=h\circ g$.
Now Theorem \ref{rational maps} tells us that there are pairwise distinct totally real algebraic numbers $\{\alpha_j\}_{j\in\mathbb{N}}$ such that
$$0\neq\widehat{h}_{g^{-1}\circ x^2 \circ g}(\alpha_j)=h(g(\alpha_j))\rightarrow 0 \quad.$$
As $g(\alpha_j)$ is in $\mathbb{Q}^{tr}(i)$ for all $j\in\mathbb{N}$, this concludes the proof.
\end{Remark}

\end{section}

\begin{section}{Chebyshev polynomials and open questions}\label{Chebyshev polynomials}

Let's go back to the quadratic polynomials $f_c = x^2 -c \in \overline{\mathbb{Q}}[x]$. We have seen in Proposition \ref{noBogomolov} that the canonical heights $\widehat{h}_{f_c}$ can get arbitrarily small on $\mathbb{Q}^{tr}$ for every rational $c\geq 2$. This behavior may change completely for non-rational $c$.

Let $q > 4$ be an element in $\mathbb{Q}\setminus\mathbb{Q}^2$. Then the Julia set of $f_{\sqrt{q}}$ is real. However, we claim that $\mathbb{Q}^{tr}$ does not have the Bogomolov property relative to $\widehat{h}_{f_{\sqrt{q}}}$.  This is due to the facts that $f_{-\sqrt{q}}$ is a Galois conjugate of $f_{\sqrt{q}}$, and that $J(f_{-\sqrt{q}})$ is not contained in the real line (see Example \ref{classification quadratic maps}). Now Theorem \ref{rational maps} proves the claim.

In the set of rational maps over $\mathbb{Q}$ with real Julia set there is one special class, namely that of the Chebyshev polynomials. 
Let $\varphi:\mathbb{C}^* \rightarrow \mathbb{C}$ be the map $x \mapsto x + x^{-1}$ and let $d$ be a natural number. We recall that the $d$-th Chebyshev polynomial is the unique polynomial $T_d$ such that the following diagram commutes.
\begin{eqnarray}\label{Chebyshev}
\begin{xy}
  \xymatrix{
      \overline{\mathbb{Q}}^* \ar[r]^{x \mapsto x^d} \ar[d]_\varphi    &   \overline{\mathbb{Q}}^* \ar[d]^\varphi  \\
      \overline{\mathbb{Q}} \ar[r]^{T_d}             &   \overline{\mathbb{Q}}   
  }
\end{xy}
\end{eqnarray}

We see at once that $f_2$ fits into this diagram for $d=2$. Hence we have $f_2 = T_2$.
An interesting fact we can deduce from Theorem \ref{rational maps} is that $\mathbb{Q}^{tr}$ has the Bogomolov property relative to the standard height $h$ coming from the map $x\mapsto x^d$, but not relative to $\widehat{h}_{T_d}$, although these heights are related in a very strong way. This relation can be easily made explicit.

\begin{Proposition}\label{Chebyshev height}
For all $z \in \overline{\mathbb{Q}}^*$ and all natural numbers $d\geq2$ we have $\widehat{h}_{T_d}(z+z^{-1})=2h(z)$.
\end{Proposition} 

\textbf{Proof:} As in \eqref{Chebyshev} we define $\varphi(x)=x+x^{-1}$. We have to check that $\frac{1}{2}\widehat{h}_{T_d}\circ\varphi$ fulfills the two conditions given in Theorem \ref{canonical height} for the canonical height $\widehat{h}_{x^d}=h$. Using the commutativity of \eqref{Chebyshev} we get $\frac{1}{2}\widehat{h}_{T_d}(\varphi(z^d))=\frac{1}{2}\widehat{h}_{T_d}(T_d (\varphi(z)))= d\frac{1}{2}\widehat{h}_{T_d}(\varphi(z))$.
As $\varphi$ has degree two, we also have $\frac{1}{2}\widehat{h}_{T_d}\circ\varphi=\frac{1}{2}h\circ\varphi + O(1)=h+O(1)$. \hfill $\square$

\begin{Definition}\rm
A Salem number is a real algebraic integer $\alpha >1$ such that all conjugates of $\alpha$ have absolute value $\leq1$ and at least one conjugate has absolute value equal to $1$. 
\end{Definition}

As one conjugate of the Salem number $\alpha$ has absolute value 1, the inverse of a conjugate is again a conjugate of $\alpha$. This implies, using the definition, that $\alpha^{-1}$ is the only real conjugate of $\alpha$ and all other conjugates lie on the unit circle. Hence $\alpha+\alpha^{-1}$ is a totally real number.

We have seen that the Bogomolov property for $\mathbb{Q}^{tr}$ does not hold relative to $\widehat{h}_{T_d}$, $d\geq 2$. The next best bound one can ask for is a bound of Lehmer strength. This means, one can ask whether there exists a positive constant $c>0$ such that $\deg(\alpha)\widehat{h}_{T_d} (\alpha)\geq c$ for all $\alpha$ in $\mathbb{Q}^{tr} \setminus \PrePer(T_d)$. This would be a quite strong result, because it would imply that the absolute value of a Salem number is bounded away from one thus proving Lehmer's conjecture for Salem numbers. This follows from Proposition \ref{Chebyshev height} and the fact that $\alpha + \alpha^{-1}$ is totally real for all Salem numbers $\alpha$. On the other hand the existence of such a bound $c$ seems to be very likely, as $\mathbb{Q}^{tr}$ has the Bogomolov property relative to all $\widehat{h}_{f_{2-\varepsilon}}$ for every algebraic $\epsilon>0$. 
\vspace{0.3cm}

Let $K$ be a number field. Recall, that $K^{ab}$ is the maximal abelian field extension of $K$. Amoroso and Zannier proved in \cite{AZ00} that $K^{ab}$ has the Bogomolov property relative to the standard height $h$. Their result also implies the Bogomolov property of these fields relative to $\widehat{h}_{T_d}$, the canonical height associated to a Chebyshev polynomial of degree at least one.

\begin{Proposition}\label{Chebyshev Bogomolov}
Let $T_d (x)$ be the $d$-th Chebyshev polynomial, where $d$ is at least 2. Let $K$ be any number field. Then the field $K^{ab}$ has the Bogomolov property relative to $\widehat{h}_{T_d}$.
\end{Proposition}

\textbf{Proof:} Let $\alpha$ be an arbitrary element in $K^{ab}\setminus \PrePer(T_d)$. Take a pre-image $\beta$ of $\alpha$ under the map $z\mapsto z+z^{-1}$. Then we have $[K^{ab}(\beta):K^{ab}]\leq 2$. From the choice of $\alpha$ we know $h(\beta)\neq 0$, hence by Proposition \ref{Chebyshev height} and \cite{AZ00}, Theorem 1.1, we get
$$\widehat{h}_{T_d}(\alpha)=2h(\beta)\geq c(K)\left(\frac{\log 4}{\log\log10}\right)^{-13}\quad,$$
for a constant $c(K)>0$ only depending on the ground field $K$. \hfill $\square$ 

\vspace{0.3cm}

The result of Amoroso and Zannier we have used here is an extension of a theorem due to Amoroso and Dvornicich. Amoroso and Dvornicich proved in \cite{AD99} that the maximal abelian field extension over $\mathbb{Q}$ has the Bogomolov property relative to $h$, answering a question raised by Zannier at a conference in Zakopane, Poland, held in honor of Schinzel's 60th birthday.

The result that $\mathbb{Q}^{ab}$ has the Bogomolov property relative to $h$ can also be stated dynamically.

\begin{center}
\textit{The field $\mathbb{Q}(\PrePer(f ))$ has the Bogomolov property relative to $\widehat{h}_{f}$, where $f=x^2$. }
\end{center}
Kronecker's theorem (Proposition \ref{Northcott} $a)$ for the map $x^d$) and Proposition \ref{Chebyshev height} show that the preperiodic points of $T_d$ are given by the set $\{\zeta + \zeta^{-1} \vert \zeta \text{ root of unity } \}$. Hence $\mathbb{Q}(\PrePer(T_d ))$ is an abelian extension of $\mathbb{Q}$ and Proposition \ref{Chebyshev Bogomolov} shows that the statement above is also true if $f$ is a Chebyshev polynomial. An interesting question is for which other rational maps $f$ this property is true.
According to the results of Habegger in \cite{Ha11} it seems to be very likely that this holds for Latt\`es maps $f$ defined over the rational numbers.

Notice that, up to linear conjugation, the maps $x^d$ and $T_d$ are the only polynomials such that infinitely many preperiodic points lie in $\mathbb{Q}^{ab}$. This was proven by Dvornicich and Zannier in \cite{DZ07}, Theorem 2.

\vspace{0.2cm}

Although we cannot prove a higher dimensional analogue of Theorem \ref{rational maps}, we will state a possible generalization as a question.

\begin{Question}
Let $(X,L,f)$ be a polarized algebraic dynamical system defined over a totally real number field $K$. Which of the following statements are equivalent?
\begin{itemize}
	\item[i)] There exists a positive constant $c$ such that $\widehat{h}_{X,L,f}$ on $X(\mathbb{Q}^{tr})$ is either zero or bounded from below by a positive constant.
	\item[ii)] There is a $\sigma \in \Gal(\overline{\mathbb{Q}}/\mathbb{Q})$ such that the Julia set $J(\sigma(f))$ is not contained in $X(\mathbb{R})$.
	\item[iii)] The set $\PrePer(f)\cap X(\mathbb{Q}^{tr})$ is not Zariski dense in $X$.
\end{itemize}
\end{Question}

\end{section}

\renewcommand{\thefootnote}{}

\footnote{Lukas Pottmeyer, Fachbereich Mathematik, TU Darmstadt, Schloßgartenstr. 7, 64289 Darmstadt, E-mail: lukas.pottmeyer@mathematik.tu-darmstadt.de}

\end{document}